\documentclass{article}

\usepackage{arxiv}

\usepackage[utf8]{inputenc} 
\usepackage[T1]{fontenc}    
\usepackage[hidelinks]{hyperref}       
\usepackage{url}            
\usepackage{booktabs}       
\usepackage{amsfonts}       
\usepackage{nicefrac}       
\usepackage{microtype}      
\usepackage{xcolor, colortbl}   
\usepackage{amssymb, amsmath}
\usepackage{float}
\usepackage{enumitem}
\usepackage{mathtools, cuted}
\usepackage{verbatim}
\usepackage{tikz,xcolor}
\usepackage{stfloats}
\usepackage{sidecap}
\usepackage{algorithm}
\usepackage[noend]{algpseudocode} 
\usepackage{array,ragged2e,multirow}
\usepackage{scrtime}
\usepackage[usestackEOL]{stackengine}
\usepackage[bb=dsserif]{mathalpha}
\usepackage{calc,array}
\usepackage{etoolbox}

\AtBeginEnvironment{tabular}{
\setlist[enumerate, 1]{wide, leftmargin=*, itemsep=0pt,label=\textbf{\roman*}, before=\vspace{-1.1\baselineskip}, after=\vspace{-1.1\baselineskip}, }
}

\newcolumntype{P}[1]{>{\RaggedRight\arraybackslash}p{#1}}

\renewcommand{\headeright}{}
\renewcommand{\undertitle}{}
\renewcommand{\shorttitle}{}

\title{Automated Radiation Therapy Patient Scheduling:  \\A Case Study at a Belgian Hospital}

\author{
  Sara Frimodig\\
  RaySearch Laboratories\\
  KTH Royal Institute of Technology\\
  Stockholm, Sweden \\
  \texttt{sarhal@kth.se} \\
  \And
 Carole Mercier \\
 Department of Radiation Oncology \\
 Iridium Netwerk \\
 Antwerp, Belgium 
   \And
 Geert De Kerf \\
 Department of Radiation Oncology \\
 Iridium Netwerk \\
 Antwerp, Belgium 
}

\begin{document}
\maketitle

\begin{abstract}
The predicted increase in the number of patients receiving radiation therapy (RT) to treat cancer calls for an optimized use of resources. To manually schedule patients on the linear accelerators delivering RT is a time-consuming and challenging task. 
Operations research (OR), a discipline in applied mathematics, uses a variety of analytical methods to improve decision-making. In this paper, we study the implementation of an OR method that automatically generates RT patient schedules at an RT center with ten linear accelerators. 
The OR method is designed to produce schedules that mimic the objectives used in the clinical scheduling while following the medical and technical constraints. The resulting schedules are clinically validated and compared to manually constructed, historical schedules for a time period of one year. It is shown that the use of OR to generate schedules decreases the average patient waiting time by 80\%, improves the consistency in treatment times between appointments by 80\%, and increases the number of treatments scheduled the machine best suited for the treatment by more than 90\% compared to the manually constructed clinical schedules, without loss of performance in other quality metrics. Furthermore, automatically creating patient schedules can save the clinic many hours of administrative work every week. 
\end{abstract}


\section{Background}

Cancer incidences are increasing globally, mainly due to a growing and aging population \cite{Ferlay2021}. 
In Europe, it is estimated that around 50\% of cancer patients would benefit from radiation therapy (RT), but that the actual use is lower \cite{Borras2015}. 
This variation largely depends on the countries' gross national income and 
results in waiting lists and delays of RT treatments \cite{CHEN2008}. 
Long waiting times for RT causes reduced effectiveness of the treatments due to tumor growth, higher risk of local recurrence, prolonged symptoms and also leads to higher levels of psychological stress among patients \cite{CHEN2008,Zumer2020,Fortin2002,Gomez2015,VanHarten2015}. Furthermore, long waiting times can also cause stress for the RT staff, which can potentially jeopardize the quality of the treatments \cite{French2004}. 
The timeliness for treatment is one of the most important quality indicators for RT \cite{VanLent2013}. Other RT quality indicators include appropriateness of care, patient experience and patient involvement, such as the possibility of choosing a time window for the treatment appointments
\cite{Harden2022, Olivotto2015}. 

An RT workflow starts with a patient consultation, followed by booking time-slots for both CT simulation and treatment delivery. CT simulation is the process where CT images are acquired for the area to be treated. The time between CT simulation and first treatment delivery is used to customize a treatment plan to each patient, taking into account  factors such as size and location of the tumor and the treatment intent. A patient-specific treatment schedule is also created.  
The majority of the patients are treated daily on an outpatient basis, meaning patients travel to the RT department on weekdays for several weeks. However, weekly and every-other-day schedules are also used clinically. 
The duration of a treatment session (\emph{fraction}) is most often 10 to 30 minutes. 
Because of stochastic patient arrivals, it is important that capacity is reserved for urgent patients that will arrive in the future.
To manually create efficient schedules for RT is both challenging and time-consuming. Therefore, the implementation of automatic scheduling algorithms can be used to support RT centers to better utilize the existing resources and improve the quality of the schedules.


Operations research (OR) is a scientific approach to improved decision-making, where complex systems and processes are analyzed and optimized using mathematical modeling, statistics, and other tools and technologies. 
OR methodologies have been widely applied to healthcare problems \cite{Rais2011,Brailsford2011,Saville2019,Vieira2016}, 
but there is not much evidence of successful implementation in clinical practice \cite{Brailsford2009, Carter2022}. 
The RT scheduling problem has been modeled using OR in several studies, for a selection see \cite{Conforti2010,Saure2012,Legrain2015,Vieira2020,Pham2021,Frimodig2022}. 
All of these models are evaluated using generated data, which is often based on clinical data, but it is not possible 
to compare the results to the clinical schedules to fully assess their suitability for clinical implementation. Only Vieira et al. \cite{Vieira2021} studies the implementation of an OR model for RT patient scheduling. Their OR model performs weekly appointment scheduling and is tested using data from two different European clinics. The model is shown to improve the quality in the weekly schedule on all key performance indicators as identified by the clinics, such as reducing time variations between appointment times by 51\%. However, it is only evaluated for a time period of one week, which is not representative to assess the dynamics in a treatment schedule, since it does not capture the uncertainty in patient arrivals or the fact that treatments span over multiple weeks.

In a previous paper, we presented a mathematical model for the RT scheduling problem \cite{Frimodig2023a}. 
Daily batch scheduling collects new patients during one day, and schedules their fractions on a machine, day and time for a planning horizon of three months. Uncertainty in future arrivals is handled by a dynamic time reservation method, which takes the expected future high-priority patients and adds them to the mathematical model as placeholder patients.

The \emph{main contribution} in this paper is the evaluation of an OR method for automatic RT patient scheduling from a clinical point of view. The potential for clinical implementation of the previously developed method \cite{Frimodig2023a} is investigated using a clinical dataset from a one-year period. The automatically generated schedules were validated by the clinical staff to ensure compliance with the medical and technical constraints posed at the clinic. Furthermore, the manually constructed clinical schedule spanning one year was compared to the automatically generated schedule for the same time period. Schedule quality metrics were developed in the clinical collaboration and were used to evaluate the schedules. We show that to automatically generate schedules would both be time-saving for the clinical staff, improve the timeliness of treatments, and increase the quality of care. 


\section{Methods}
\label{Sec:clinical_collab}
The potential for clinical implementation of the mathematical model for RT patient scheduling presented by Frimodig et al. \cite{Frimodig2023a} is evaluated using clinical input from different stakeholders. The input from medical physicists, booking administrators and radiation oncologists is used for data acquisition, defining quality metrics, and schedule validation. 

\subsection{The RT center: Iridium Netwerk}
\label{Sec:Iridium}
Iridium Netwerk is a large RT center with four different hospitals located in Antwerp, Belgium. In 2020, ten linacs were in operation, and approximately 5500 RT treatments were delivered to 4000 patients.
The main location (S1) has four linacs, and the three satellite locations (S2-S4) operate two linacs each. In 2020, three different linac types were operational. 
Two linacs of the same type are called \textit{beam-matched}; both machines can deliver the same treatment plan, and a machine switch between fractions is possible. If the matched machines are at the same hospital, they are \textit{completely} matched, otherwise the match is \textit{partial}. 
A switch between partially matched machines is allowed only in special cases.
The yearly schedule at Iridium Netwerk includes weekends and holidays, where no treatments are delivered, and days for planned maintenance or other quality assurance activities. When a machine is planned to be unavailable, the affected fractions can either be scheduled on a beam-matched machine the same day or postponed. For the latter, a fraction can be added at the end of the treatment or the patient can be treated twice on a day if the time between two fractions is at least 6 hours. 
In 2020, 34\% of all weekdays had some type of machine unavailability planned due to holidays, planned maintenance or other quality assurance activities.

Each patient is appointed to a \textit{treatment protocol} that contains the scheduling instructions, such as required linac type, the duration of each fraction, dose prescription, and the minimum number of fractions to be treated each week. 
Since time is needed to prepare the treatment, the protocol defines the minimum number of days between CT simulation and treatment start. 
The most common treatment protocols at Iridium in 2020 can be seen in Table~\ref{tab:protocols}. 
In case a patient needs two different \emph{consecutive treatments}, a secondary plan should ideally start the day after the primary plan has ended. These treatments are most common for breast cancer patients when the primary plan is followed by a boost plan. 
If the secondary plan has a different treatment protocol than the primary, it can also have different machine requirements and durations, and must therefore be handled separately in the scheduling. Approximately 17\% of all patients at Iridium Netwerk in 2020 were treated with two or more consecutive treatments.
\begin{table}
\caption{The most common treatment protocols Iridium Netwerk}
\setlength\extrarowheight{1pt}
\label{tab:protocols}
\begin{tabular}{P{1.7cm} P{1.05cm} P{2.1cm}  P{1.97cm} P{4.5cm} P{2.66cm}}
\hline
Protocol & Priority & Minimum fractions/week & Pre-treatment phase (days) & Preferred machines &Allowed machines\\
\hline
Palliative & A & 1 & 0 & M1, M2, M3, M4, M5, M6, M7, M8 & M10\\
Breast & C & 3 & 7 & M1, M4, M5, M6, M8 & M2, M3, M7, M10 \\
Prostate & C & 3  & 9 & M1, M3, M4, M5, M6, M7, M8 & M2, M10 \\
Head-Neck  & A &  5 & 11 & M2, M3, M5, M6, M10 & M1, M4 \\
Lung & B & 4  & 9  & M2, M3, M5, M6, M7, M10 & M1, M4, M8 \\
\hline
\end{tabular}
\end{table}

Each treatment protocol is associated with a priority as seen in Table~\ref{tab:protocols}. The purpose of the priority is to estimate the urgency for treatment; priority~A-patients should start treatment as soon as possible, whereas priority~B-patients can wait a few days and priority~C-patients a few weeks. In 2020 at Iridium, 37\% of the patients were priority~A, 16\% were priority~B, and 46\% were priority~C.  Since the treatments often span multiple weeks, there is a need to ensure that capacity is reserved for priority~A-patients that will arrive in the future. At Iridium Netwerk, this is handled by the booking administrator, who reserves empty timeslots for urgent patients on each machine based on historical data.

When scheduling treatment appointments, patient wishes are taken into consideration, such as preferred location and time. The booking administrators have estimated that 80\% of the patients
have a preference for treatment time, of which 65\% prefer a treatment before noon and 35\% prefer the afternoon. 
Each patient is notified of their treatment schedule soon after admission and all fractions are communicated at once and will only be re-planned in case of an unforeseen event such as a machine failure. The booking administrators estimate that it takes approximately one booking administrator's full time job to schedule the arriving patients on all sites.

\subsection{Data Acquisition}
\label{Sec:data_acquisition}
Information of all treatments in 2020 has been extracted from Iridium Netwerk's oncology information system (OIS) ARIA \cite{Varian}. The patient data has been anonymized to comply with the General Data Protection Regulation (GDPR). 

\begin{figure}[b]
\vspace{-0.3cm}
\center
\includegraphics[scale=0.7]{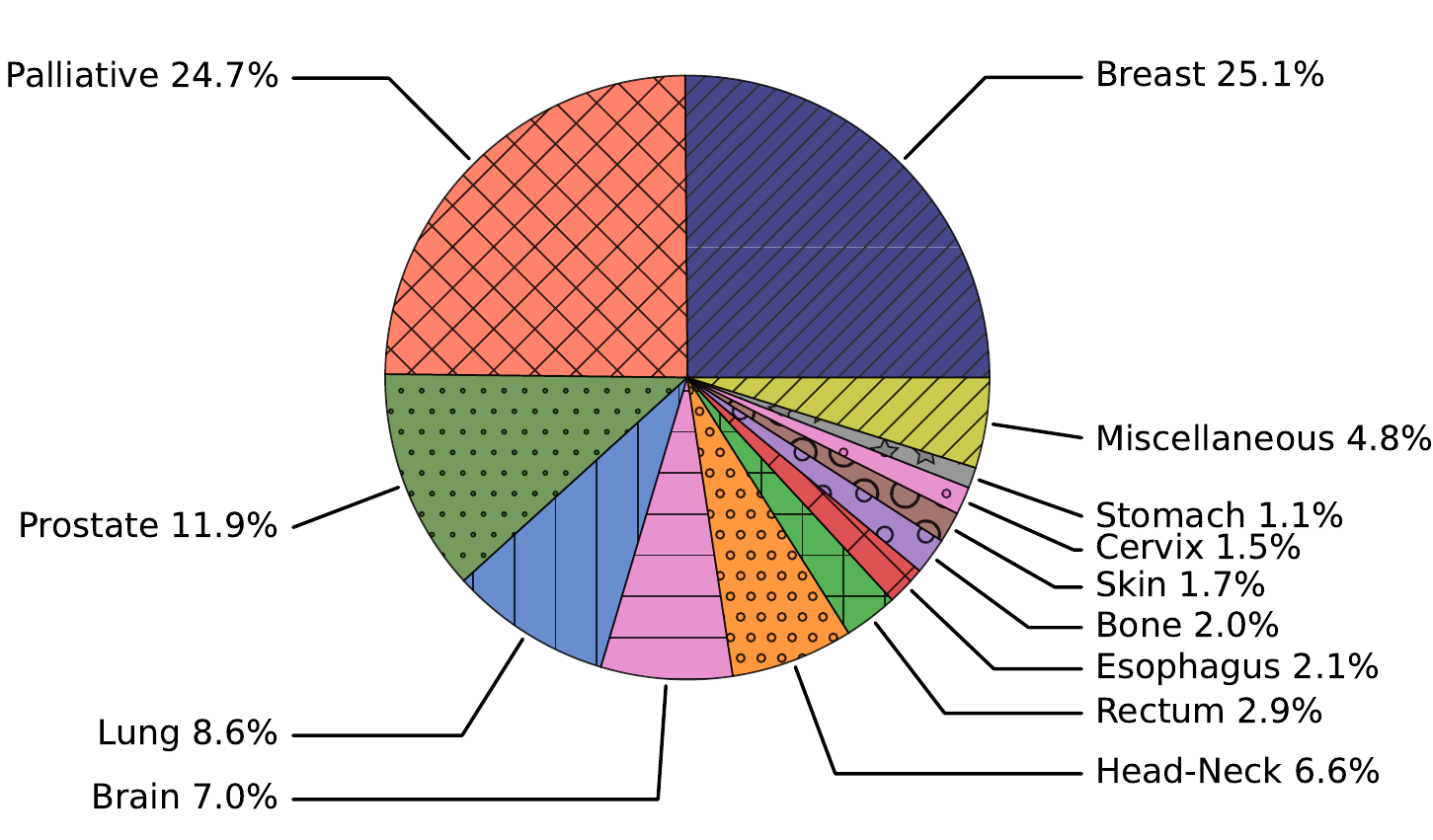}
\caption{Breakdown of cancer category at Iridium Netwerk in 2020}
\vspace{-10pt}
\label{Fig:pie_chart}
\end{figure}

The information about what \textit{treatment protocol} each patient has followed is recorded in multiple free text fields in the data. 
The free text often directly states the treatment protocol, and otherwise the prescribed dose and what machine the treatment has been delivered on are used for the protocol mapping.
The palliative patients are grouped as one category. A breakdown of the cancer categories at Iridium Netwerk in 2020 can be seen in Figure~\ref{Fig:pie_chart}. Each treatment protocol is assigned a priority by the medical staff at Iridium Netwerk.
Furthermore, each treatment is assigned a \textit{CourseID} based on PatientID, treatment protocol and the times for the delivered fractions. This is used to separate the treatments for patients that were treated with multiple treatment courses in 2020. 
When each treatment course has a CourseID, it is possible to link together \textit{consecutive treatments}. This is done based on the delivery of the fractions: if a second CourseID starts within a few days after a another CourseID for the same patient, it is assumed that they are connected. 
All patients that arrived in 2019, but had treatments scheduled in 2020, are seen as fixed in the \textit{input schedule} for the 2020 scheduling. 

The patients' \textit{hospital site preferences} are communicated verbally and not recorded in the data. According to the Iridium staff, the majority of patients want to be treated on the hospital closest to where they live. To keep the data anonymous, instead of using the patient's home address, a script assigns a site preference to each patient based on the patient's zip code and measured distance to each hospital.

The dates for  \textit{planned unavailability} in the yearly schedule 2020 are given as input to the algorithm. Moreover, the extracted clinical data has information about \textit{machine failures} and other unplanned interruptions, and the durations of these interruptions are computed and given as input at the day of the interruption. 

\begin{figure}
\center
 \includegraphics[scale=0.85]{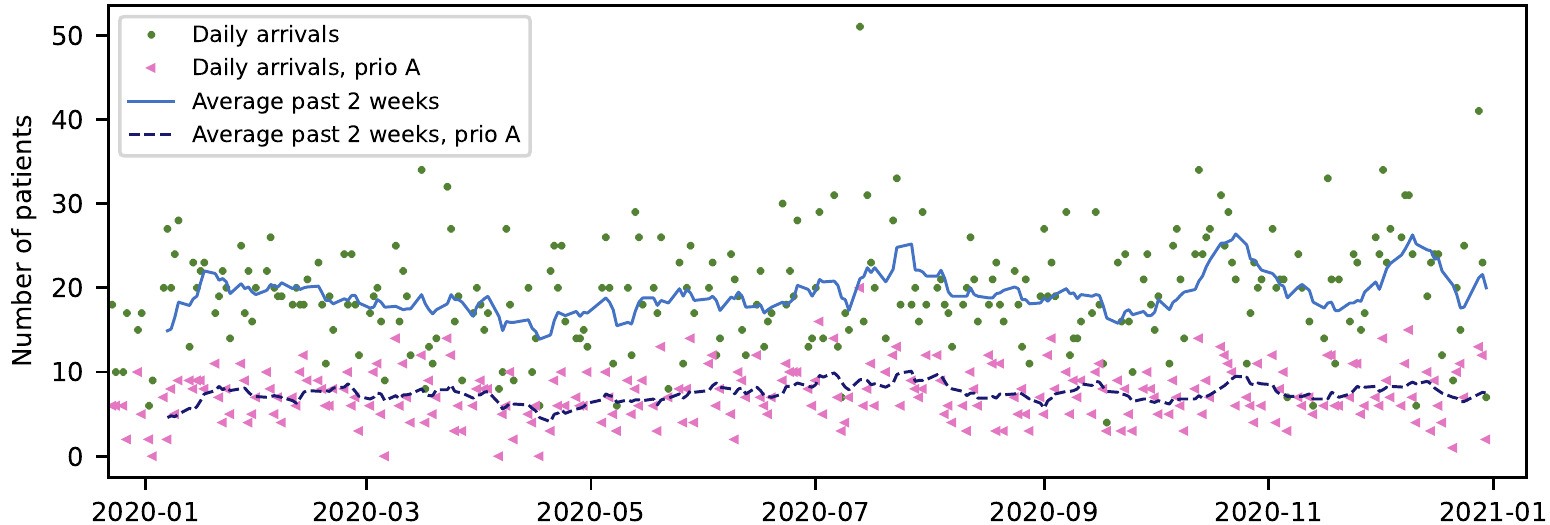}
\caption{Daily arrivals in 2020 at Iridium Netwerk, excluding weekends and holidays}
\label{Fig:Arrivals}
\end{figure}

The automatic algorithm reads information about the \textit{patient arrivals}, where each treatment has a PatientID, CourseID, creation date, treatment protocol, number of fractions, the length of the first and subsequent fractions, and site preference. For secondary treatments, it is stated what CourseID it follows. The creation date in the OIS data is assumed to mimic the arrival to the hospital. Table~\ref{tab:patient_arrivals} shows examples of patient arrivals. Figure~\ref{Fig:Arrivals} shows the daily arrivals of patients in 2020, where the rolling average during the past two weeks varies between 13.9 to 26.4 patients per day, and the average priority A patient arrivals vary between 4.1 to 10.1 patients per day. The day with the most arrivals in the year was the 13\textsuperscript{th} of July, with 51 new arrivals, of which 20 were priority A. Some examples of the input data for three weeks in 2020 can be seen in Table~\ref{tab:data_inputs}. The clinical data used by the automatic scheduling algorithm is publicly available\footnote{Access through this link: \url{https://osf.io/j2bxp/?view_only=e1402382b67f4ad0a4b8a3f4ed28088a}}.

\begin{table}
\caption{The first lines in the input data-file of patient arrivals after pre-processing}
\setlength\extrarowheight{2pt}
\label{tab:patient_arrivals}
\begin{tabular}{P{0.08\linewidth} P{0.082\linewidth} P{0.08\linewidth} P{0.175\linewidth} P{0.075\linewidth}  P{0.075\linewidth} P{0.075\linewidth}P{0.045\linewidth} P{0.083\linewidth}}
\hline
PatientID &  CourseID & Creation Date & Protocol & Number fractions & Duration \#1 (min) & Duration (min)  & Site pref & Follows CourseID \\
\hline
00001 & 11730 & 20-01-02& Bladder VMAT & 30 & 24 & 12 & S2 & -  \\
00002 & 11755 & 20-01-02& Head-Neck VMAT& 35 & 24 & 12 & S4 & -  \\
00003 & 16282 & 20-01-02& Prostate VMAT & 20 & 24 & 12 & S2 & -  \\
00004 & 14140 & 20-01-02& Liver SBRT & 8 & 50 & 50 & S3 & -  \\
00005 & 18671 & 20-01-02& Prostate VMAT & 20 & 24 & 12 & S2 & -  \\
00006 & 12402 & 20-01-03& Breast tang. fields	 & 15 & 24 & 12 & S4 & - \\
00006 & 15930 & 20-01-03& Breast photon boost & 5 & 20 & 12 & S4 & 12402  \\
\hline
\end{tabular}
\end{table}

\begin{table}[ht]
\caption{Data for the patients arriving during three random weeks, shown as example of model input}
\setlength\extrarowheight{2pt}
\label{tab:data_inputs}
\begin{tabular}{P{0.357\linewidth} P{0.18\linewidth} P{0.18\linewidth}P{0.18\linewidth}}
\hline
      Week                                &  24-28 Feb 2020 & 4-8 May 2020 & 13-17 Oct 2020 \\
\hline
Number of patient arrivals	  & 80 & 69 & 114 \\
Average session duration, minutes (SD) & 14.5 (5.8) & 15.4 (6.5) & 15.3 (7.5)\\
Average number of fractions (SD)  & 18.8 (9.1) & 19.9 (11.2) & 19.9 (10.5) \\
Percentage with consecutive treatments & 18.8\% & 20.3\%  & 17.5\%\\
Patient site preference distribution & S1: 27\% & S1: 26\% & S1: 18\%\\
                                                              & S2: 25\% & S2: 22\% & S2: 25\% \\
                                                              & S3: 29 \%& S3: 36\% & S3: 33\%\\
                                                              & S4: 19 \%& S4: 16\% & S4: 24\%  \\  
Most common treatment protocols & Breast: 26.3\% & Breast: 29.0\% & Palliative: 23.7\% \\
                                                                & Palliative: 23.7\% & Palliative: 27.5\%& Breast: 23.7\%\\
                                                                & Prostate: 16.2 \% & Prostate: 13.0\% & Prostate: 14.9\% \\
\hline
\end{tabular}
\end{table}

\subsection{Automatic Scheduling Approach}
\label{Sec:auto_setup}

The mathematical model developed in \cite{Frimodig2023a} is based on a column generation approach. It generates a schedule where each treatment session is assigned to a machine, day and time window, and the patient is assigned a specific start time within the time window in a post-processing step. The automatic scheduling is done at the end of each day taking patients from previous days into account. 
Daily batch scheduling can leverage accumulated information of the patients that have arrived during the day to reduce patient waiting time and overdue time \cite{Pham2021}. 

At Iridium Netwerk each patient is notified of their treatment schedule soon after admission. However, literature shows that the majority of patients find it reasonable to receive a notification of their start date three days in advance \cite{Olivotto2015}. In \cite{Frimodig2023a} a compromise is used; for priority~B and C patients, the notification period is one week and priority~A patients are notified immediately. Iridium Netwerk communicates all fractions at once, so the same approach is applied to the automatic algorithm. A patient schedule can change until is has been communicated to the patient, i.e., booking decisions are postponed to the next day if patients are scheduled after the notification period.

The automatic algorithm reserves machine capacity for high priority patients to account for uncertainty in future patient arrivals. In \cite{Frimodig2023a}, a dynamic time blocking method is developed, in which the expected value of the future patient arrivals are included in the model as placeholder patients. This allows for a trade-off between the current (actual) patients to be scheduled, and the urgent patients that are expected to arrive in the future. Experiments in \cite{Frimodig2023a} show that this method better utilizes the resources than using a static time reservation method, which is used in the manual scheduling.

When a machine fails it results in unexpected downtime. 
At Iridium Netwerk, the booking administrators reschedule affected patients on beam-matched machines if possible, or otherwise move the fractions to other days. 
In the automatic scheduling, it is assumed that the daily interruptions are known already in the morning. The machine failures are managed in a pre-processing step by a heuristic algorithm that blocks failed machines for new bookings. The affected patients are moved to other machines or days according to what is allowed in the treatment protocol, and the list of the current patients is thereafter scheduled.

\subsection{Schedule Quality Metrics}
\label{Sec:QualityMetrics}
Quality indicators have been published for certain tumor types, including minimizing waiting times or completing treatment within a certain time frame \cite{Leroy2019}. The patient experience is mostly affected by treatment-related information, professional standard and short time-interval between diagnosis and treatment \cite{Petersen2015}. It is also affected by the available treatment options, such as the possibility of choosing a time window for the treatment appointments \cite{Olivotto2015}. 
Since different cancer centers have different clinical setups, the quality metrics can differ between different clinics. 
To evaluate the automatically generated schedules and compare them to the manually constructed ones, multiple schedule quality metrics have been developed in collaboration with Iridium Netwerk through interviews with booking administrators, medical physicists, radiation oncologists and managers. Table~\ref{tab:quality_metrics} shows the quality metrics relevant at Iridium Netwerk. Each quality metric is presented as an objective that the booking administrators aim to fulfill, together with a short description and an associated priority that should reflect the importance of the objective. For example, the booking administrators always try to schedule the treatments on one of the best suited machines (objective \ref{enum:iv}), however, if the waiting time (objective \ref{enum:i}) is shorter on another machine, that is always preferred as that objective has a higher priority. 
\begin{table}[ht]
\small
\caption{Schedule quality metrics and their priorities}  
\setlength\extrarowheight{2pt}
\label{tab:quality_metrics}
\begin{tabular}{P{0.39\linewidth} P{0.47\linewidth} P{0.06\linewidth}} 
\hline
Objective &  Description & Priority \\
\hline
\begin{enumerate}
\item \label{enum:i} Minimize waiting time between referral and treatment start
\end{enumerate} & The waiting time should be minimized for all patients, but is most important for high priority patients & 1 \\
\begin{enumerate}
\setcounter{enumi}{1}
\item \label{enum:ii} Minimize the time deviations between fractions for each patient 
\end{enumerate}  & Consistency of appointments & 4 \\
\begin{enumerate}
\setcounter{enumi}{2}
\item \label{enum:iii} Minimize the number of fractions scheduled outside of preferred time windows 
\end{enumerate} & For patients that have preferred treatment times & 4 \\
\begin{enumerate}
\setcounter{enumi}{3}
\item  \label{enum:iv} Minimize the number of fractions scheduled on non-preferred machines 
\end{enumerate} & Among the list of allowed machines for each treatment, some are preferred over the others & 3 \\
\begin{enumerate}
\setcounter{enumi}{4}
\item \label{enum:v} Minimize the number of switches between partially beam-matched machines  
\end{enumerate} & Number of switches between same-type machines that are at different hospital sites should be minimized & 3\\
\begin{enumerate}
\setcounter{enumi}{5}
\item \label{enum:vi} Minimize fractions scheduled outside of patient's preferred hospital site 
\end{enumerate} & Fulfill patients' site preferences regarding treating hospital & 2\\
\begin{enumerate}
\setcounter{enumi}{6}
\item \label{enum:vii} Minimize overall treatment time 
\end{enumerate}& Gap days in the schedule can be induced when there is machine unavailability, but should be kept to a minimum & 1 \\
\hline
\end{tabular}
\end{table}

\section{Results}
\label{Sec:results}
In a simulation of year 2020, the actual patient arrivals from Iridium Netwerk are used as input to the mathematical algorithm. Since daily batch scheduling is performed, there is one schedule to create each of the 254 workdays in 2020. The planning horizon for each day is three months. The patients that have arrived during the day are scheduled together with a list of unscheduled patients from the days before, while considering the input schedule that describes the partial occupation on each machine due to previous scheduling decisions. The experiments are run on a computer with Windows 10 that has an Intel\textsuperscript{\textregistered} Core{\texttrademark} i9-7940X X-series processor and 64~GB of RAM. The mathematical model for automatic schedule generation is created in Python~3.8 and solved using IBM ILOG CPLEX 20.1 in the Python API.

To evaluate the potential for clinical implementation of the automatically generated schedules, three factors must be assessed; the \emph{validity} of the schedules, 
the \emph{quality} of the schedules, 
and the \emph{time} it takes to construct the schedules. 

\subsection{Schedule Validation}
\label{Sec:model_validation}
To ensure the automatically generated schedules can be implemented in practice, two different types of validation have been done by the staff at Iridium Netwerk. The  \emph{patient schedule validation} checks that the treatment is scheduled correctly for each patient, and the \emph{machine schedule validation} controls that the schedule is valid on each machine. Since it would be too time consuming to manually confirm the fulfillment of all medical and technical constraints for all ten machines and all 4000 treated patients at Iridium Netwerk in 2020, the schedule validation is done for a subset of patients and machines. First, a set of 15 randomly selected patient schedules are validated, covering a range of different treatment protocols. The patient information is reviewed together with the resource information data to confirm that the schedule would work in practice. The reviewed patient information includes fractionation schemes, number of fractions, fraction duration, time for pre-treatment and minimum number of fractions per week. Second, a random machine is selected for a time period of two weeks. The inspected machine information includes double bookings, machine unavailability, machine validity for each patient's protocol and that the machine capacity is not exceeded.

In addition to the manual validation performed by the medical staff at Iridium Netwerk, an automatic framework for validating the schedules has also been developed. It tests all the automatically generated schedules, both machine schedules (capacity, planned and unplanned unavailability), and patient schedules. For the patient schedules, every medical and technical constraint is tested for every patient and their treatment protocols; number of fractions, allowed start days, maximum fractions per day, minimum fractions per week, maximum interval between fractions, all fractions scheduled on beam-matched machines, all fractions scheduled on allowed machine, minimum number of days for pre-treatment, fractionation schemes, non-conventional protocol fulfillment, fractions at least 6 hours apart if two scheduled on same day, and finally the time interval between primary and secondary treatment for consecutive treatments.

\subsection{Schedule Quality}
\label{Sec:NumResults}
\begin{figure}
\center
\includegraphics[scale=0.85]{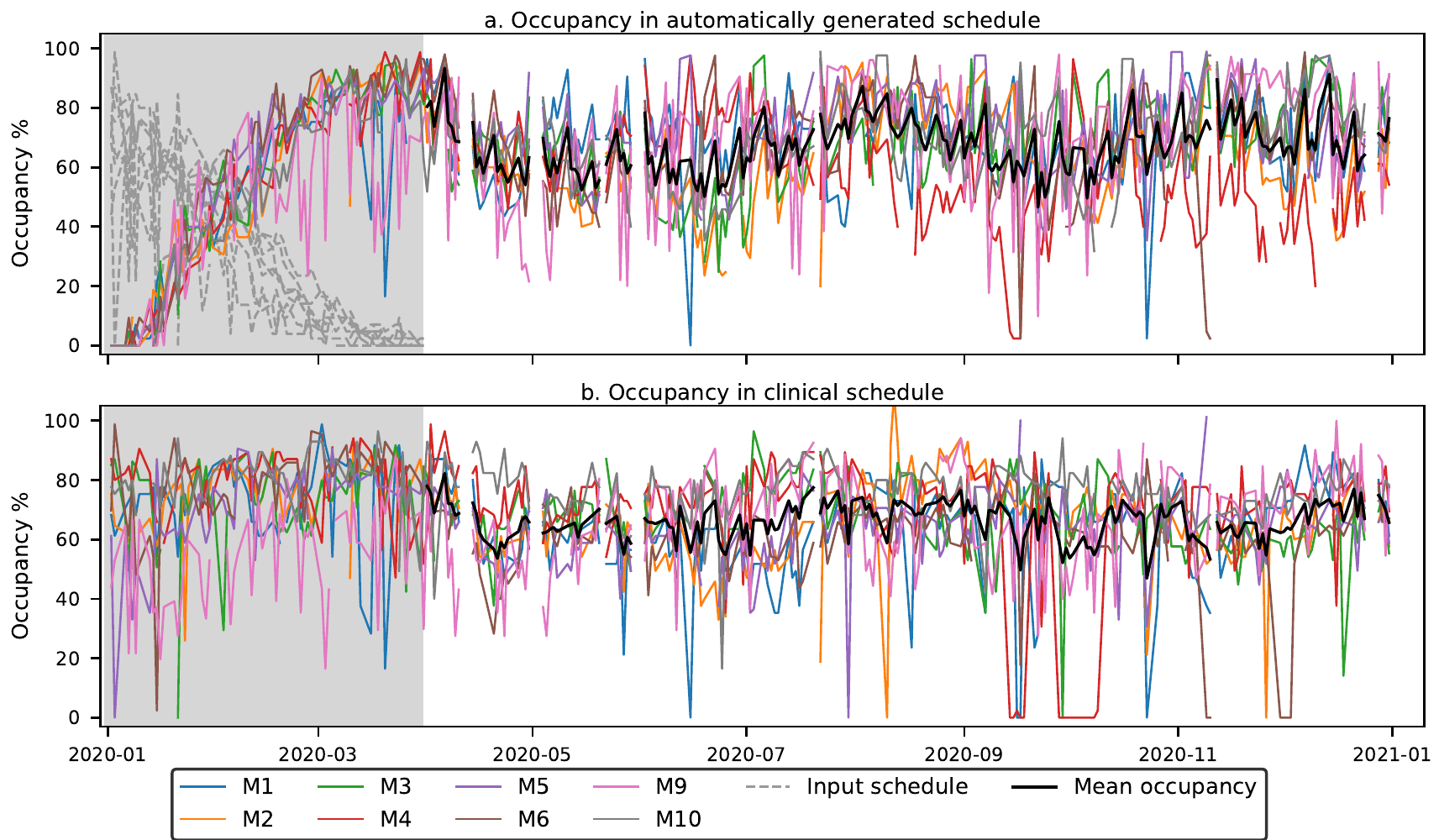}
\caption{Occupancy in the automatically generated schedule and in the clinical schedule. The gray region is the first three months, where the input schedule with manually scheduled patients from 2019 make up a proportion of the total occupancy in the optimized schedule. The days where no data is shown are unavailable machine dates}
\label{Fig:Manual_opt_sched_occupancy}
\end{figure}
Because a shortage of staff at one of the hospital sites at Iridium Netwerk, the machines at that site were not fully operational in 2020. Therefore, machines M7 and M8, together with all patients that were scheduled on those machines in the clinical schedule, are not included in the comparison between the schedules.

All patients that arrived to Iridium Netwerk in 2019, but had treatments scheduled in 2020, are seen as fixed in the input schedule for the automatic scheduling in 2020. These treatments span approximately the first three months in 2020, and therefore this period is excluded when comparing the manual and automatic schedules, since it is in fact a mix of the two. 
Figure \ref{Fig:Manual_opt_sched_occupancy} shows the occupancy (i.e., the machine utilization rate) on each machine for the automatically generated schedule and the actual clinical schedule in 2020. Figure~\ref{Fig:Manual_opt_sched_occupancy}a also shows the occupancy resulting from the 2019 manually scheduled patients in gray scale. 
From Figure \ref{Fig:Manual_opt_sched_occupancy}, it can be seen that the average overall occupancy on the machines is similar between the automatic and the manual scheduling: for all machines, it is 66.9\% for the automatic and 65.9\% for the manual schedule in 2020, excluding the first three months. 

\begin{figure}
\center
\includegraphics[scale=0.9]{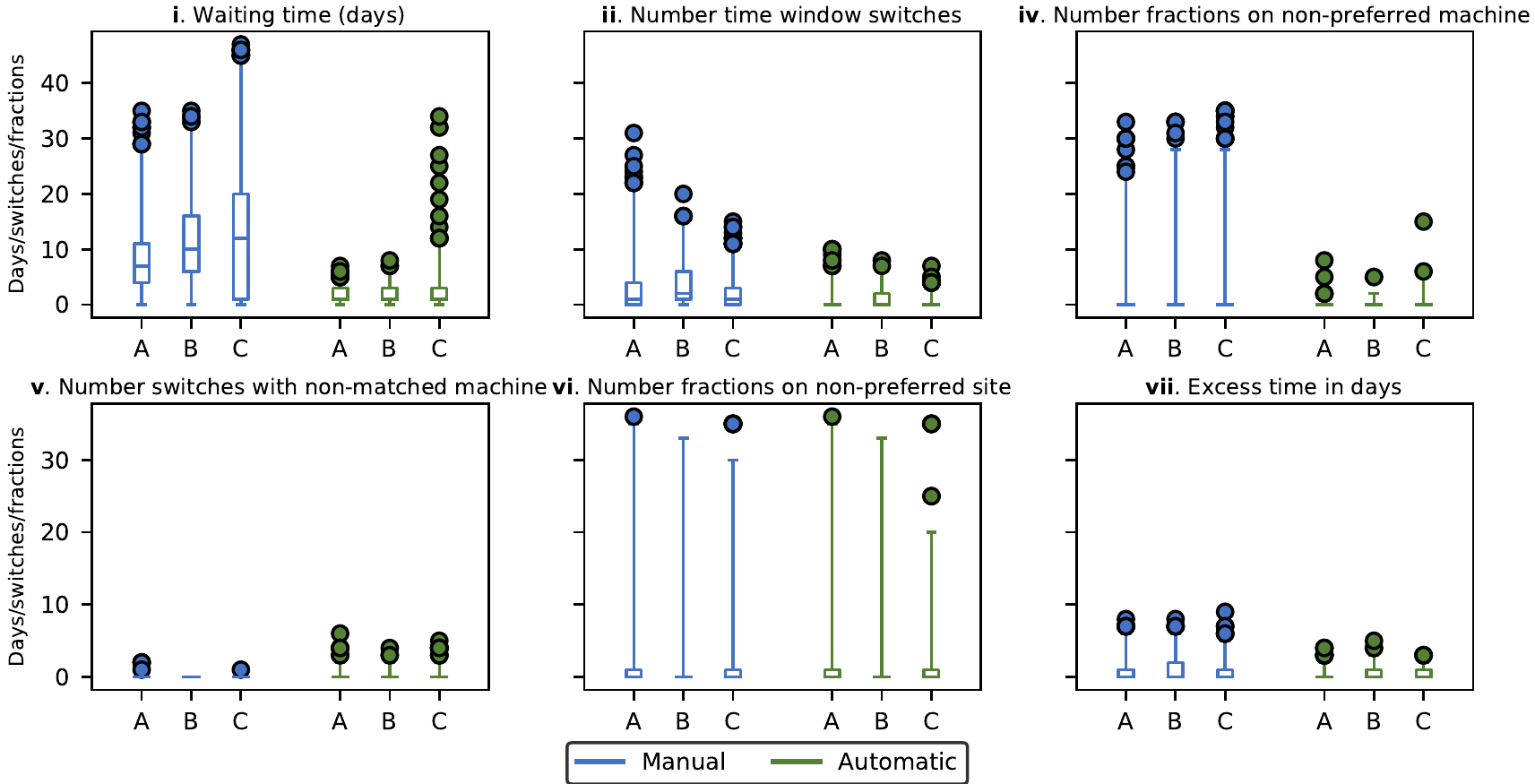}
\caption{Comparison between costs for priority group A, B and C for automatically created schedules and manually constructed clinical schedules, with top 1\% marked as outliers}
\label{Fig:results_manual_automatic}
\end{figure}

To assess the quality of the automatically generated schedules, the schedule quality metrics presented in Section~\ref{Sec:QualityMetrics} were evaluated and compared to the historical schedules at Iridium Netwerk. When analyzing the manual schedules, a number of patients were identified to have unusually long waiting times. These patients were reviewed by the staff at Iridium Netwerk, and 169 of the patients had notes about why the treatment was delayed, usually because of a request from the patient, or due to clinical need, or because of reasons related to COVID-19. Therefore, these patients were excluded from the quality metric comparison altogether. Figure~\ref{Fig:results_manual_automatic} shows the results for measures that correspond to six of the quality metrics for the manually constructed clinical schedules and the automatically generated schedules. 
The reason why objective \ref{enum:iii}, to minimize the number of fractions scheduled outside the patient's preferred time window, is not included in Figure~\ref{Fig:results_manual_automatic} is because Iridium Netwerk did not record the patients' time window preferences, thus it is not possible to measure what method that best fulfilled this objective.

In Figure \ref{Fig:results_manual_automatic}, similar performance is observed for objective \ref{enum:vi}, to minimize the number of fractions scheduled on a non-preferred hospital site. 
For objectives \ref{enum:ii}, \ref{enum:iv} and \ref{enum:vii}, Figure \ref{Fig:results_manual_automatic} shows that the automatic algorithm fulfills the goals better than the clinical schedules, both if including the 1\% outliers or when excluding them. The number of time window switches for a patient (objective \ref{enum:ii}) was on average 2.8 in the clinical schedule and 0.5 in the automatically generated schedule, which means that the consistency in appointment times was 84\% better in the optimized schedule. The average number of fractions scheduled on a non-preferred machine (objective \ref{enum:iv}) decreased by 93\% between the manual and automatic schedule (from 1.3 to 0.1 on average per patient). The excess time in days, which corresponds to the number of extra gap days in the schedule (objective \ref{enum:vii}) decreased from an average of 0.7 in the clinical schedule to 0.3 in the optimized schedule, a 55\% decline. 
Objective \ref{enum:v}, to minimize the number of switches between machines that are only partially beam-matched, is the only objective where the manual schedule outperforms the automatically generated one; the average number of switches to partially beam-matched machines is 0.003 in the manual schedule, and 0.14 is the optimized schedule, which is still very low. 

Finally, Figure \ref{Fig:results_manual_automatic} illustrates that for objective \ref{enum:i}, to minimize the waiting times, the automatically generated schedule outperforms the manually constructed schedule for all priority groups, but especially for urgent patients. The waiting time is computed as the time from the earliest allowed start date, taking into account the required number of days of the pre-treatment phase. For priority group~A, the average waiting time in the clinical schedule was 8.1 days, compared to 1.5 days  in the automatic scheduling approach, which corresponds to a 81\% decrease in average waiting time. For priority~B, the decrease in average waiting time was 86\% (from 11.3 to 1.6 days), and for priority~C-patients the average waiting time decreased 86\% (from 13.3 to 1.8 days). Furthermore, the maximum waiting times decreased from 35 to 7 days for priority~A-patients, from 35 to 8 for priority~B and from 47 to 34 for patients in priority group C.

\section{Discussion}
It has been shown that the previously developed OR method for the RT scheduling problem can produce clinically acceptable schedules. The quality of the generated schedules is better than the historical schedules based on the quality metrics that have been developed within the clinical collaboration. Since the waiting times are shown to be decreased by the automatic algorithm, the time between first consultation and end of treatment is reduced by more than 80\% while the length of the pre-treatment phase remains the same. Three main reasons why the automatic scheduling algorithm produces schedules that are superior to the clinical schedules have been identified. First, the automatic algorithm performs daily batch scheduling instead of scheduling the patients one-by-one. This offers more degrees of freedom since the algorithm can leverage patient information accumulated throughout the day, and it has also been shown in previous studies that this leads to better schedules \cite{Pham2021}. Secondly, the handling of uncertainty in future urgent patients is improved in the automatic algorithm; it uses a dynamic time reservation method that allows for a tradeoff between actual patients and the patients that are expected to arrive in the future, which gives shorter waiting times than the static time reservation method that is used clinically today \cite{Frimodig2023a}. Third, the automatic algorithm uses mathematical optimization, and is designed to maximize the performance in the schedule quality metrics described in Section~\ref{Sec:QualityMetrics}. It is shown in \cite{Frimodig2023a} that, if the problem is modeled correctly, the resulting schedule will be close to optimal using the proposed mathematical approach. Another potential reason why the automatic schedule has shorter waiting times than the clinical schedule lies in the result of objective \ref{enum:v}; the automatic algorithm sometimes schedules a patient with planned switches between hospitals if it results in shorter waiting times. 


As previously stated, approximately 170 outlier patients with unusually long waiting times were excluded from the quality measure plot in Figure~\ref{Fig:results_manual_automatic}. However, it is possible that additional waiting times are long in the clinical schedule not as a result of the scheduling strategy, but as a result of shortage of staff or other resources, due to personal requests, or influenced by the COVID pandemic (either because patients were sick or because they chose to delay their treatments). If excluding the top 15\% waiting times in the clinical schedule, the automatically generated schedule still decreases the average waiting times by 75\% for priority~A (6.2 to 1.5 days), 82\% for priority~B (8.8 to 1.6 days), and 81\% for priority~C (9.5 to 1.8 days). Furthermore, due to the shortage of staff in one of the hospital sites, the two machines on that site together with the patients scheduled on them in the clinical schedule were not included when creating the automatic schedules or in the quality comparison. If including these patients and machines with full capacity, the results are very similar, thus, this does not seem to be a large source of error.

For each treatment course, the creation date in the OIS is assumed to mimic the arrival to the hospital. This will be correct in the majority of the cases since the patient is entered into the OIS at arrival. However, there can be cases when it is not necessarily decided to treat with RT already at arrival (for example patients first treated with surgery or chemotherapy, where it is sometimes decided later to also treat with RT). There is no way to know this from the available data, but the Iridium staff assure that these cases are very uncommon. Another source of error is that in the automatic scheduling, it is assumed that the daily interruptions are known already in the morning, even though it is likely that the failure is not known until it actually happens. However, for the rescheduling of the affected patients this does not make a difference, since we limit the rescheduled fractions to be after the failure occurs.

Since the hospital site preferences were not recorded, they were instead estimated from the patient zip codes. 
This is a source of error when evaluating the algorithm for this quality metric, but should not be a problem if used in practice when the actual preferences are known. Moreover, different cancer centers can have different scheduling objectives. 
In some countries, there are nationally adopted waiting time targets stating the maximum waiting time between referral and start of treatment. In a previous paper, Frimodig et al. \cite{Frimodig2022} showed that there is a close correlation between minimizing the waiting time and minimizing the violations of the waiting time targets, thus, if this objective is relevant it is straightforward to include in the mathematical model.

At Iridium Netwerk, approximately one full time job is spent every week on scheduling patients on the linacs. For each daily batch scheduling, the automatic algorithm is allowed to run for maximum one hour, 
and the average computation time is 33 minutes. Since the scheduling is done once daily, the automatic algorithm will spend at most five hours per week constructing schedules, a process that can be done in the background. 
Although the model is designed to capture all the major medical and technical constraints present in the clinical scheduling, some patients may have patient-specific needs that are not implemented in the model. Therefore, the generated schedule may need to be manually adjusted. 
Furthermore, the staff may also need to spend time validating the schedules. Even with time spent on these actitivies, there is a potential to save many hours of administrative work every week by automatically generating schedules. 

From a clinical point of view, the medical staff at Iridium Netwerk expresses a large interest in applying the automatic scheduling algorithm clinically. They especially see two main benefits of applying the model, where the first is that the clinical schedules are optimized, and the second is that the algorithm can be used to simulate department changes, such as buying a new linac, which can help make more informed decisions. Another main advantage is the rescheduling of patients in case of machine breakdowns, where an algorithm can take into account all the different scheduling parameters (overall treatment time, minimum number of fractions per week, etc), something that is difficult for a human to do efficiently.
However, several challenges exist for actual clinical implementation. The most important is that a user interface is needed, preferably  developed within the OIS, or at least with communication to it. That way, the information of planned maintenance and protocol data already exist within the system. The OIS should also hold the information of the treatment protocol for each patient, removing the need for the data acquisition steps. Furthermore, to enable all the schedule quality metrics to be incorporated in the automatic algorithm, the patients' preferences for time windows and hospital sites should be recorded. In the current solution, the CPLEX software is used for the mathematical optimization. Since this is an expensive system, the possibility of switching to an open source solver should be investigated. 

\section{Conclusions}
\label{Sec:conclusions}
To manually schedule patients for RT is complex and labor-intensive. As the number of cancer cases in the world are rapidly increasing, automated scheduling algorithms can  optimize the resources to achieve short waiting times. In this paper, a previously developed operations research method for automated RT scheduling is evaluated from a clinical point of view in a collaboration with Iridium Netwerk.  
The automatically generated schedules are clinically validated and compared to the historical schedules from Iridium Netwerk for a time period of one year. 
The results show that when compared to the clinical schedule, the automatically generated schedule decreases the average patient waiting time by approximately 80\%, increases the time consistency between treatment appointments by 80\%, increases the number of fractions scheduled on preferred machines by 90\%, and also improves the majority of the other schedule quality metrics substantially. 




\bibliographystyle{unsrt}  
\bibliography{mybibliography}

\begin{thebibliography}{10}

\bibitem{Ferlay2021}
Jacques Ferlay, Murielle Colombet, Isabelle Soerjomataram, Donald~M Parkin,
  Marion Pi{\~{n}}eros, Ariana Znaor, and Freddie Bray.
\newblock {Cancer statistics for the year 2020: An overview}.
\newblock {\em International Journal of Cancer}, 149(4):778--789, aug 2021.

\bibitem{Borras2015}
Josep~M. Borras, Yolande Lievens, Peter Dunscombe, Mary Coffey, Julian Malicki,
  Julieta Corral, Chiara Gasparotto, Noemie Defourny, Michael Barton, Rob
  Verhoeven, Liesbeth {Van Eycken}, Maja Primic-Zakelj, Maciej Trojanowski,
  Primoz Strojan, and Cai Grau.
\newblock {The optimal utilization proportion of external beam radiotherapy in
  European countries: An ESTRO-HERO analysis}.
\newblock {\em Radiotherapy and Oncology}, 116(1):38--44, 2015.

\bibitem{CHEN2008}
Zheng Chen, Will King, Robert Pearcey, Marc Kerba, and William~J. Mackillop.
\newblock The relationship between waiting time for radiotherapy and clinical
  outcomes: A systematic review of the literature.
\newblock {\em Radiotherapy and Oncology}, 87(1):3 -- 16, 2008.

\bibitem{Zumer2020}
Barbara {\v{Z}}umer, Maja {Pohar Perme}, Simona Jereb, and Primo{\v{z}}
  Strojan.
\newblock {Impact of delays in radiotherapy of head and neck cancer on
  outcome}.
\newblock {\em Radiation Oncology}, 15(1):202, dec 2020.

\bibitem{Fortin2002}
Andr{\'{e}} Fortin, Isabelle Bairati, Michele Albert, Lynne Moore, Jos{\'{e}}e
  Allard, and Christian Couture.
\newblock Effect of treatment delay on outcome of patients with early-stage
  head-and-neck carcinoma receiving radical radiotherapy.
\newblock {\em International Journal of Radiation Oncology Biology Physics},
  52(4):929--936, 2002.

\bibitem{Gomez2015}
Daniel~R Gomez, Kai-Ping Liao, Stephen~G Swisher, George~R Blumenschein,
  Jeremy~J. Erasmus, Thomas~A Buchholz, Sharon~H Giordano, and Benjamin~D
  Smith.
\newblock Time to treatment as a quality metric in lung cancer: Staging
  studies, time to treatment, and patient survival.
\newblock {\em Radiotherapy and Oncology}, 115(2):257--263, 2015.

\bibitem{VanHarten2015}
Michel~C {Van Harten}, Frank~J.P. Hoebers, Kenneth~W Kross, Erik~D {Van
  Werkhoven}, Michiel~W.M. {Van Den Brekel}, and Boukje~A.C. {Van Dijk}.
\newblock Determinants of treatment waiting times for head and neck cancer in
  the netherlands and their relation to survival.
\newblock {\em Oral Oncology}, 51(3):272--278, 2015.

\bibitem{French2004}
Helen~C. French.
\newblock {Occupational stresses and coping mechanisms of therapy radiographers
  – a qualitative approach}.
\newblock {\em Journal of Radiotherapy in Practice}, 4(1):13--24, jun 2004.

\bibitem{VanLent2013}
W.~A.M. {Van Lent}, R.~D. {De Beer}, B.~{Van Triest}, and W.~H. {Van Harten}.
\newblock {Selecting indicators for international benchmarking of radiotherapy
  centres}.
\newblock {\em Journal of Radiotherapy in Practice}, 12(1):26--38, 2013.

\bibitem{Harden2022}
Susan~V. Harden, Kim‐Lin Chiew, Jeremy Millar, and Shalini~K. Vinod.
\newblock {Quality indicators for radiation oncology}.
\newblock {\em Journal of Medical Imaging and Radiation Oncology},
  66(2):249--257, mar 2022.

\bibitem{Olivotto2015}
I.A. Olivotto, J.~Soo, R.A. Olson, L.~Rowe, J.~French, B.~Jensen, A.~Pastuch,
  R.~Halperin, and P.T. Truong.
\newblock {Patient Preferences for Timing and Access to Radiation Therapy}.
\newblock {\em Current Oncology}, 22(4):279--286, aug 2015.

\bibitem{Rais2011}
Abdur Rais and Ana Viana.
\newblock {Operations Research in Healthcare: a survey}.
\newblock {\em International Transactions in Operational Research},
  18(1):1--31, jan 2011.

\bibitem{Brailsford2011}
Sally Brailsford and Jan Vissers.
\newblock {OR in healthcare: A European perspective}.
\newblock {\em European Journal of Operational Research}, 212(2):223--234, jul
  2011.

\bibitem{Saville2019}
Christina~E. Saville, Honora~K. Smith, and Katarzyna Bijak.
\newblock {Operational research techniques applied throughout cancer care
  services: a review}.
\newblock {\em Health Systems}, 8(1):52--73, 2019.

\bibitem{Vieira2016}
Bruno Vieira, Erwin~W Hans, Corine {Van Vliet-Vroegindeweij}, Jeroen {Van De
  Kamer}, and Wim {Van Harten}.
\newblock Operations research for resource planning and -use in radiotherapy: a
  literature review.
\newblock {\em BMC Medical Informatics and Decision Making}, 16(149), 2016.

\bibitem{Brailsford2009}
S.~C. Brailsford, P.~R. Harper, and M.~Pitt.
\newblock {An analysis of the academic literature on simulation and modelling
  in health care}.
\newblock {\em Journal of Simulation}, 3(3):130--140, 2009.

\bibitem{Carter2022}
Michael~W. Carter and Carolyn~R. Busby.
\newblock {How can operational research make a real difference in healthcare?
  Challenges of implementation}.
\newblock {\em European Journal of Operational Research}, apr 2022.

\bibitem{Conforti2010}
D.~Conforti, F.~Guerriero, and R.~Guido.
\newblock {Non-block scheduling with priority for radiotherapy treatments}.
\newblock {\em European Journal of Operational Research}, 201(1):289--296,
  2010.

\bibitem{Saure2012}
Antoine Saur\'e, Jonathan Patrick, Scott Tyldesley, and Martin~L. Puterman.
\newblock Dynamic multi-appointment patient scheduling for radiation therapy.
\newblock {\em European Journal of Operational Research}, 223(2):573 -- 584,
  2012.

\bibitem{Legrain2015}
Antoine Legrain, Marie-Andr{\'{e}}e Fortin, Nadia Lahrichi, and Louis-Martin
  Rousseau.
\newblock Online stochastic optimization of radiotherapy patient scheduling.
\newblock {\em Health Care Management Sciences}, 18:110--123, 2015.

\bibitem{Vieira2020}
Bruno Vieira, Derya Demirtas, Jeroen~B. van~de Kamer, Erwin~W. Hans,
  Louis-Martin Rousseau, Nadia Lahrichi, and Wim~H. van Harten.
\newblock {Radiotherapy treatment scheduling considering time window
  preferences}.
\newblock {\em Health Care Management Science}, 23(4):520--534, dec 2020.

\bibitem{Pham2021}
Tu-San Pham, Louis-Martin Rousseau, and Patrick {De Causmaecker}.
\newblock {A two-phase approach for the Radiotherapy Scheduling Problem}.
\newblock {\em Health Care Management Science}, 25(2):191--207, jun 2022.

\bibitem{Frimodig2022}
Sara Frimodig, Per Enqvist, Mats Carlsson, and Carole Mercier.
\newblock Comparing optimization methods for radiation therapy patient
  scheduling using different objectives.
\newblock {\em arXiv:2211.01150}, 2022.

\bibitem{Vieira2021}
Bruno Vieira, Derya Demirtas, Jeroen~B. van~de Kamer, Erwin~W. Hans, Willem
  Jongste, and Wim van Harten.
\newblock {Radiotherapy treatment scheduling: Implementing operations research
  into clinical practice}.
\newblock {\em PLoS ONE}, 16(2 February):1--13, 2021.

\bibitem{Frimodig2023a}
Sara Frimodig, Per Enqvist, and Jan Kronqvist.
\newblock A column generation approach for radiation therapy patient scheduling
  with planned machine unavailability and uncertain future arrivals.
\newblock {\em arXiv:2303.10985}, 2023.

\bibitem{Varian}
{Varian Medical Systems, Inc.}
\newblock
  \url{https://www.varian.com/products/software/information-systems/aria-ois-medical-oncology}.
\newblock 2020.

\bibitem{Leroy2019}
Roos Leroy, Cindy {De Gendt}, Sabine Stordeur, Geert Silversmit, Leen Verleye,
  Viki Schillemans, Isabelle Savoye, Katrijn Vanschoenbeek, Joan Vlayen,
  Liesbet {Van Eycken}, Claire Beguin, C{\'{e}}cile Dubois, Laurens Carp, Jan
  Casselman, Jean-Fran{\c{c}}ois Daisne, Philippe Deron, Marc Hamoir, Esther
  Hauben, Olivier Lenssen, Sandra Nuyts, Carl {Van Laer}, Jan Vermorken, and
  Vincent Gr{\'{e}}goire.
\newblock {Quality indicators for the management of head and neck squamous cell
  carcinoma}.
\newblock {\em KCE Reports}, 305, 2019.

\bibitem{Petersen2015}
Gitte~Stentebjerg Petersen, Janne~Lehmann Knudsen, and Mette~Marianne Vinter.
\newblock {Cancer patients' preferences of care within hospitals: A systematic
  literature review}.
\newblock {\em International Journal for Quality in Health Care},
  27(5):384--395, 2015.

\end{thebibliography}

\end{document}